\documentclass[12pt,a4paper, reqno]{amsart}
\usepackage[utf8]{inputenc}
\usepackage[T1]{fontenc}
\usepackage{amsmath}
\usepackage{amsthm}
\usepackage{amssymb}
\usepackage[abbrev]{amsrefs}
\usepackage{mathrsfs}
\usepackage[dvipsnames]{xcolor}
\usepackage{bm}
\usepackage{enumitem}
\usepackage{hyperref}
\usepackage{graphicx}
\usepackage[all]{xy}
\AtBeginDocument{\def\MR#1{}}
\makeatletter
\@namedef{subjclassname@2020}{%
\textup{2020} Mathematics Subject Classification}
\makeatother
\numberwithin{equation}{section}
\allowdisplaybreaks
\newtheorem{thmintro}{}

\newtheorem{Theoremintro}[thmintro]{Theorem}
\newtheorem{Theorem}{Theorem}[section]

\newtheorem{Proposition}[Theorem]{Proposition}
\theoremstyle{remark}

\newtheorem{Question}[Theorem]{Question}

\newcommand{\xx}{\ensuremath{\bm{x}}}
\newcommand{\XX}{\ensuremath{{X}}}
\newcommand{\YY}{\ensuremath{Y}}
\newcommand{\Mt}{\ensuremath{\mathcal{M}}}
\newcommand{\Nt}{\ensuremath{\mathcal{N}}}
\newcommand{\Pt}{\ensuremath{\mathcal{P}}}
\newcommand{\RR}{\ensuremath{\mathbb{R}}}
\newcommand{\NN}{\ensuremath{\mathbb{N}}}
\newcommand{\VV}{\ensuremath{V}}
\newcommand{\F}{\ensuremath{\mathcal{F}}}
\newcommand{\XB}{\ensuremath{\mathcal{X}}}
\newcommand{\LL}{\ensuremath{\mathscr{L}}}
\newcommand{\Id}{\ensuremath{\mathrm{Id}}}
\newcommand{\Lip}{\operatorname{\rm{Lip}}}
\newcommand{\abs}[1]{\left\lvert#1\right\rvert}
\newcommand{\norm}[1]{\left\lVert#1\right\rVert}
\newcommand{\enangle}[1]{\left\langle#1\right\rangle}

\newcommand{\enpar}[1]{\left(#1\right)}

\begin{document}
\title[Envelope maps and duals of Lipschitz-free $p$-spaces]{On Banach envelopes and duals\\ of Lipschitz-free $\bm{p}$-spaces for $\bm{0<p<1}$}
\author[F. Albiac]{Fernando Albiac}
\address{Department of Mathematics, Statistics, and Computer Sciences--Ina\-Mat$^2$ \\
Universidad P\'ublica de Navarra\\
Campus de Arrosad\'{i}a\\
Pamplona\\
31006 Spain}
\email{fernando.albiac@unavarra.es}
\author[J. L. Ansorena]{Jos\'e L. Ansorena}
\address{Department of Mathematics and Computer Sciences\\
Universidad de La Rioja\\
Logro\~no\\
26004 Spain}
\email{joseluis.ansorena@unirioja.es}
\begin{abstract}
With the aim to better understand the intricate geometry of the class of Lipschitz free $p$-spaces $\F_p(\Mt)$ when $0<p<1$, in this note we study their Banach envelopes and prove that if $0<p<1$ and $\Mt$ is a metric space then the Banach envelope map of $\F_p(\Mt)$ is one-to-one, thus solving in the positive a problem raised by Kalton in \cite{AlbiacKalton2009}. This property has important applications to the linear structure of this family of spaces, being the most immediate one that the dual space of $\F_p(\Mt)$ separates the points of $\F_p(\Mt)$.
\end{abstract}
\thanks{The authors acknowledge the support of the Spanish Ministry for Science and Innovation under Grant PID2022-138342NB-I00 for \emph{Functional Analysis Techniques in Approximation Theory and Applications}}
\maketitle
\section{Introduction and background}\noindent
The Banach envelope is a fundamental construction that allows a quasi-Banach space to be embedded into a Banach space in a universal way, preserving its linear structure while endowing it with the geometric regularity of local convexity. Banach envelopes were explicitly introduced in the 1970s by Peetre in \cite{Peetre1974}*{p.\@ 125} and Shapiro \cite{Shapiro1977}*{p.\@ 116}. Since then, they have been considered by many authors in different spaces and settings. Without intending to be exhaustive, we refer the reader to \cites{MN2002,CRS2007,KaltonSukochev2008,KamLin2014,CW2017} for some recent advances on the subject.

Recall that the \emph{Banach envelope} of a quasi-Banach space $X$ consists of a Banach space $\widehat{X}$ together with a linear contraction $J_X\colon X \to \widehat{X}$ satisfying the following property: for every Banach space $Y$ and every continuous linear map $T\colon X \to Y$ there is a unique continuous linear map $\widehat{T}\colon \widehat{X}\to Y$ such that $\widehat{T}\circ J_X=T$, that is,
\[
\xymatrix{
\widehat{X}\ar[drr]^{\widehat{T}} & \\
X \ar[u]^{J_X} \ar[rr]_T && Y
}
\]
and the ``extension'' $\widehat{T}$ has a norm bounded by the norm of $T$.

In particular, $X$ and $\widehat{X}$ have the same dual space. We will refer to $J_{X}$ as the \emph{(Banach) envelope map} of $X$. Since the Banach envelope of $X$ is defined by means of a universal property, it is unique in the sense that if a Banach space $Z$ and a bounded linear map $J\colon X\to Z$ satisfy the property, then there is an linear isometry $S\colon \widehat{X} \to Z$ with $S\circ J_X=J$.

Roughly speaking, the map $J_{X}$ determines which parts of the original space survive convexification and how much of its geometry and functional structure is retained. In general, the envelope map need not be injective. When it fails to be so, the quasi-Banach space contains nonzero elements that are annihilated in every Banach context, rendering them invisible to bounded linear functionals and obstructing any meaningful transfer of analytic structure. Conversely, injectivity ensures that the space embeds faithfully into its Banach envelope, making it possible to import tools from Banach space theory into the study of quasi-Banach spaces.

It is well known (see \cite{AABW2021}*{\S 9} or \cite{AAW2021c}*{Lemma 2}) that the envelope map $J_{X}$ is injective if and only if the dual space $X^{\ast}$ separates the points of $X$. Assuming the point-separation property is essential in some situations. Indeed, finding sufficient conditions that guarantee that the envelope map is one-to-one has both theoretical and practical interest. To substantiate our claim, let us explain a situation that occurs in applications. Some quasi-Banach spaces $X$ (such as Lipschitz free $p$-spaces for $0<p<1$) are naturally constructed from a non-complete quasi-normed space $X_0$ that is dense in $X$. This definition leads to a precise description on how $X_0$ embeds into $\widehat{X}$ via $J_X$, but not on how $X$ does! Thus, if we want to know that the character of the members of $X_{0}$ inside $\widehat{X}$ is preserved in passing to its completion, it is important to know whether we can transfer properties from $J_X|_{X_0}$ to $J_X$. In particular, the following question arises.

\begin{Question}\label{PhilQ}
Let $X$ be a quasi-Banach space and let $X_{0}$ be a dense subspace of $X$. Assume that $X^{*}$ separates the points of $X_{0}$. Does it follow that $X^{*}$ separates the points of $X$? Equivalently, assuming that $J_X|_{X_0}\colon X_{0}\to \widehat{X}$ is one-to-one, is $J_{X}\colon X\to \widehat{X}$ one-to-one?
\end{Question}

The authors were glad to learn about other colleagues concern with this kind of problem. For instance, in his study of the Banach envelope of $\ell_{1,\infty}$ (\cite{Pietsch2009}), Pietsch raised and solved the following question, which he elevated to the category of `philosophical’.

\begin{Question}[\cite{Pietsch2009}*{p.\@ 214}]
Let $X_0$ and $Z$ be normed spaces. Suppose that $Z$ is complete and that $T\colon X_0\to Z$ is a one-to-one linear operator. Is the extension of $T$ to the completion of $X_0$ one-to-one?
\end{Question}

In turn, a negative answer to Question~\ref{PhilQ} was provided in \cite{AAW2021c}*{Theorem 2}, where it was shown that if $X$ is a separable quasi-Banach space with $N=(X^*)^\perp\not=\{0\}$ and $X/N$ infinite-dimensional, then there exists a dense subspace $X_{0}$ of $X$ such that $J_{X}$ is one-to-one on $X_{0}$ but fails to be one-to-one on $X$.

These considerations on the injectivity of the envelope map are particularly relevant in emerging nonlinear settings, such as Lipschitz-free $p$-spaces over metric (or $p$-metric) spaces for $0<p<1$. Given a pointed $p$-metric space $\Mt$ there is a unique $p$-Banach space $\F_p(\Mt)$ such that $\Mt$ embeds isometrically in $\F_p(\Mt)$ via a canonical map denoted $\delta_{p,\Mt}$, and for every $p$-Banach space $\XX$ and every Lipschitz map $f\colon \Mt\to\XX$ with Lipschitz constant $\Lip(f)$ that maps the base point in $\Mt$ to the zero vector in $\XX$ extends to a unique bounded linear map $E[f;p]\colon \F_p(\Mt)\to\XX$ with $\Vert E[f;p]\Vert =\Lip(f)$. Pictorially we have
\begin{equation}\label{Dibujo1}
\xymatrixcolsep{5pc}\xymatrix{
\F_p(\Mt) \ar[rd]^{E[f;p]}&\\
\Mt \ar[u]^{\delta_{p,\Mt}} \ar[r]_{f}& \XX.
}
\end{equation}

The class of Lipschitz-free $p$-spaces $\F_p(\Mt)$ for $0<p<1$ was introduced in \cite{AlbiacKalton2009}, where they were used to give the first-known examples for each $0<p<1$ of separable $p$-Banach spaces which are Lipschitz-isomorphic but not linearly isomorphic. Subsequently, a systematic study of the geometry of the spaces $\F_p(\Mt)$ was initiated in \cite{AACD2018}, where some foundational problems posed in \cite{AlbiacKalton2009} were solved and new intriguing questions were raised. Let us introduce some terminology to provide a context for what follows.

Let $0$ denote the base point of $\Mt$. Let $\RR_{0}^{\Mt}$ be the space of all maps $f\colon \Mt\to\RR$ so that $f(0)=0$, and let $\Pt(\Mt)$ be the linear span in the algebraic dual $(\RR_{0}^{\Mt})^{\#}$ of the evaluations $\delta_{\Mt}(x)$, where $x$ runs through $\Mt$, defined by
\begin{equation}\label{Delta}
\enangle{\delta_{\Mt}(x),f}=f(x), \quad f\in \RR_{0}^{\Mt}.
\end{equation}
Note that $\delta_{\Mt}(0)=0.$

If $\mu=\sum_{x\in \Mt} a_x\delta_\Mt(x)\in\Pt(\Mt)$, put
\begin{equation}\label{definition}
\norm{\mu}_{\F_p(\Mt)}=\sup\norm{\sum_{x\in\Mt}a_x f(x)}_\XX,
\end{equation}
the supremum being taken over all $p$-normed spaces $\enpar{\XX,\norm{\cdot}_\XX}$ and all $1$-Lipschitz maps $f\colon \Mt\to\XX$ with $f(0)=0$. It is straightforward to check that formula \eqref{definition} defines a $p$-seminorm on $\Pt(\Mt)$. In fact, more can be said, which answers an early question from \cite{AlbiacKalton2009} on Lipschitz free $p$-spaces.

\begin{Proposition}[\cite{AACD2018}*{Proposition 4.1 and Theorem 4.10}]\label{thm:AACD}
Let $(\Mt,d)$ be a $p$-metric space, $0<p\le 1$. Then $\norm{\cdot}_{\F_p(\Mt)}$ is a $p$-norm over $\Pt(\Mt)$. Besides, for all $\mu\in\Pt(\Mt)$,
\[
\norm{\mu}_{\F_p(\Mt)}=\inf\enpar{\sum_{j\in F}\abs{a_j}^p}^{1/p},
\]
the infimum being taken over all finite families $(a_j)_{j\in F}$ in $\RR$ for which there are $(x_j)_{j\in F}$ and $(y_j)_{j\in F}$ in $\Mt$ such that $x_j\not=y_j$ for all $j\in F$, and
\begin{equation}\label{eq:expansion}
\mu=\sum_{j\in F} a_j\frac{\delta_\Mt(x_j)-\delta_\Mt(y_j)}{d(x_j,y_j)}.
\end{equation}
\end{Proposition}

The Lipschitz free $p$-space $\F_p(\Mt)$ is the completion of the $p$-normed space $(\Pt(\Mt), \norm{\cdot}_{\F_p(\Mt)})$, and the canonical embedding $\delta_{p,\Mt}$ of $\Mt$ into $\F_p(\Mt)$ is just $\delta_\Mt$ regarded as a map into $\F_p(\Mt)\supseteq \Pt(\Mt)$. In the case when $p=1$, the spaces $\F_p(\Mt)$ become the Lipschitz-free spaces over $\Mt$, also known, depending on tastes and schools, as Arens--Eells spaces or transportation cost spaces, and denoted by $\F(\Mt)$, which have been widely studied from several angles (see e.g. \cites{WeaverBook2018,GodefroyKalton2003,OO19}).

It is known (see \cite{AACD2018}*{Proposition 4.20}) that if $\Mt$ is a metric space, then the Banach envelope of $\F_p(\Mt)$ is $(\F(\Mt),J_{p}[\Mt])$, where $J_{p}[\Mt]$ is just the operator $E[\delta_{1,\Mt};p]$ that we obtain when we consider the following particular case of the diagram in \eqref{Dibujo1} with $f=\delta_{1,\Mt}$:
\[
\xymatrixcolsep{5pc}\xymatrix{
\F_p(\Mt) \ar[rd]^{J_p[\Mt]}&\\
\Mt \ar[r]_{\delta_{1,\Mt}}\ar[u]^{\delta_{p,\Mt}}& \mathcal F(\Mt).
}
\]

The envelope map $J_{p}[\Mt]$ is injective on a dense subspace of $\F_p(\Mt)$ by Proposition~\ref{thm:AACD}. However, this does not guarantee that $J_{p}[\Mt]$ is one-to-one on the entire space $\F_p(\Mt)$ as we already discussed. Kalton was well aware of this drawback in the advancement of the theory of Lipschitz free $p$-spaces and made it explicit in \cite{AlbiacKalton2009}*{Section 4}. This article is devoted to answering that question in the positive by proving the following theorem.

\begin{Theoremintro}[Main Theorem]\label{thm:LFp1Envelope}
Let $\Mt$ be a metric space and $0<p<1$. Then the envelope map $J_{p}[\Mt]\colon\F_p(\Mt)\to\F(\Mt)$ is one-to-one.
\end{Theoremintro}

The consequences of this result are far-reaching. It guarantees that classical tools from Banach space theory can be meaningfully applied to $\F_p(\Mt)$ through its envelope, and it ensures the separation of points by continuous linear functionals.
\section{Preparatory results}\noindent
Our proof of Theorem~\ref{thm:LFp1Envelope} takes advantage of the recent advances on another key question about Lipschitz free $p$-spaces. To state it, we recall that these spaces provide a canonical linearization process of Lipschitz maps between $p$-metric spaces. Indeed, given two $p$-metric spaces $\Mt$ and $\Nt$ and a Lipschitz map $h\colon\Nt\to\Mt$ with $h(0)=0$, the bounded linear map
\[
L[h;p]:=E[\delta_{p,\Mt}\circ h;p]
\]
makes the diagram
\[
\xymatrixcolsep{3pc}\xymatrix{
\F_p(\Nt) \ar[r]^{L[h;p]} &\F_p(\Mt)\\
\Nt \ar[r]_{h} \ar[u]^{\delta_{p,\Nt}} &\Mt \ar[u]_{\delta_{p,\Mt}}
}
\]
commute. Note that $\norm{L[h;p]}=\Lip(h)$.

Although significant progress has been made in the understanding of Lipschitz-free $p$-spaces (see \cites{AACD2019,AACD2021,AACD2022}), the possibility to linearize isometric embeddings between these spaces when $p<1$ is an unexploited tool due to the fact that its validity is not ensured.

Given a pointed $p$-metric space $\Mt$ and $0\in\Nt \subset\Mt$, we denote by
\[
L[\Nt,\Mt;p]\colon\F_p(\Nt)\to\F_p(\Mt)
\]
the canonical linearization of the inclusion map of $\Nt$ into $\Mt$. In this terminology, we write the question mentioned above as follows.

\begin{Question}[\cite{AACD2018}*{Question 6.2}]\label{qt:LIEP}
Let $0<p<1$, $\Mt$ be a pointed $p$-metric space, and $0\in\Nt \subset\Mt$. Is $L:=L[\Nt,\Mt;p]$ and isomorphic embedding? If the answer is positive, is there a constant $C=C(p)$ independent of $\Nt$ and $\Mt$ such that $\norm{ L^{-1}}\le C$?
\end{Question}

Cuth and Raunig \cite{CuthRaunig2024} gave the following partial answer to Question~\ref{qt:LIEP}.
\begin{Theorem}[see \cite{CuthRaunig2024}*{Theorem 4.21}]\label{thm:CR}
Let $0<p\le 1$, $\Mt$ be a metric space, and $0\in\Nt \subset\Mt$. Then, $L:=L[\Nt,\Mt;p]$ is an isomorphic embedding, and $\norm{L^{-1}}\le 1500\cdot 18^{1/p}$.
\end{Theorem}

We will also take advantage of some important features of Lipschitz-free $p$-spaces over Euclidean spaces obtained in \cite{AACD2022}. Specifically, we will use the following structural result.

\begin{Theorem}[see \cite{AACD2022b}*{Theorem 3.9}]\label{thm:FpVSchauder}
Let $0<p\le 1$ and $d\in\NN$. Then $\F_p(\RR^d)$ has a Schauder basis.
\end{Theorem}

A \emph{Markusevich basis} of a quasi-Banach space $\XX$ is a family $\XB=(\xx_\gamma)_{\gamma\in\Gamma}$ in $\XX$ for which there is $\XB^*=(\xx_\gamma^*)_{\gamma\in\Gamma}$ in $\XX^*$ such that
\begin{itemize}
\item $\xx_\gamma^*(\xx_\alpha)=\delta_{\alpha,\gamma}$ for all $\alpha$, $\gamma\in\Gamma$,
\item the linear span of $\XB$ is dense in $\XX$, and
\item the linear span of $\XB^*$ is weak$^*$-dense in $\XX^*$, that is, the coefficient transform defined on $\XX$ by
\[
f\mapsto \enpar{\xx_\gamma^*(f)}_{\gamma\in\Gamma},
\]
is one-to-one.
\end{itemize}
Quasi-Banach spaces with a Markusevich basis clearly have the point-separation property. Since Schauder bases are in particular Markusevich bases, we have the following.

\begin{Proposition}\label{prop:SchauderSeparates}
Let $\XX$ be a quasi-Banach space with a Schauder basis. Then the envelope map of $\XX$ is one-to-one.
\end{Proposition}

We point out that Proposition~\ref{prop:SchauderSeparates} is the most commonly used tool to show the injectivity of the envelope map. However, proving that a Lipschitz free $p$-space over a given separable metric space has a Schauder basis has proven to be a challenging task. Thus, Proposition~\ref{prop:SchauderSeparates} would be of help to show Theorem~\ref{thm:LFp1Envelope} only when combined with other techniques.

The third ingredient we will use in the proof of Theorem~\ref{thm:LFp1Envelope} is the theory of $\LL_\infty$-spaces. Recall that a Banach space $\XX$ is an \emph{$\LL_\infty$-space} if and only if there is a constant $C$ such that for every finite-dimensional subspace $\VV\subset\XX$ there are $d\in\NN$ and linear maps $P\colon\XX\to \ell_\infty^d$ and $J\colon \ell_\infty^d \to \XX$ such that $\norm{J} \norm{P} \le C$, $P\circ J=\Id_{\ell_\infty^d}$, and $\VV\subset J(\ell_\infty^d)$. If $D\in[1,\infty)$ is such that we can choose any $C>D$, we say that $\XX$ is an $\LL_{\infty,D}$-space.

By the principle of small perturbations, $\XX$ is an $\LL_{\infty,D}$-space if and only if there is a directed set $(P_\lambda,J_\lambda,d_\lambda)_{\lambda\in \Lambda}$ such that for each $\lambda\in \Lambda$ and $d_\lambda\in\NN$, $P_\lambda\colon\XX\to \ell_\infty^{d_\lambda}$ and $J\colon \ell_\infty^{d_\lambda} \to \XX$ are linear operators such that $P\circ J=\Id_{\ell_\infty^{d_\lambda}}$, $\limsup_{\lambda\in \Lambda} \norm{J_\lambda} \norm{P_\lambda} \le \lambda$, $(J_\lambda(\ell_\infty^{d_\lambda}))_{\lambda\in \Lambda}$ is increasing, and $\cup_{\lambda\in\Lambda} J_\lambda(\ell_\infty^{d_\lambda})$ is dense in $\XX$. Hence, given an arbitrary set $I$, $c_0(I)$ is an $\LL_{\infty,1}$-space. By the principle of local reflexivity, its bidual space $\ell_\infty(I)$ is an $\LL_{\infty,1}$-space as well (see \cites{LinPel1968,LinRos1969}).

Any metric space $\Mt$ isometrically embeds in a suitable Banach space $\XX$ and, actually, we can choose $\XX=\F(\Mt)$. In turn, $\XX$ linearly and isometrically embeds in $\ell_\infty(B_{\XX^*})$ by the Hahn--Banach theorem. For further reference, we record the well-known consequence of these observations.

\begin{Theorem}\label{thm:metricembedding}
Any metric space isometrically embeds into an $\LL_{\infty,1}$-space.
\end{Theorem}
\section{Proof of Main Theorem}\noindent
We are now ready to complete the proof of Theorem~\ref{thm:LFp1Envelope}. Prior to it, we introduce some additional terminology.

If $\Mt$ is a quasi-metric space, $\Nt$ is a subspace of $\Mt$, and $\mu\in\Pt(\Mt)$ admits an expansion as in \eqref{eq:expansion} with $x_j$, $y_j\in \Nt$ for all $j\in F$, we say that $\mu$ is \emph{supported} on $\Nt$.

\begin{proof}[Completion of the Proof of Theorem~\ref{thm:LFp1Envelope}]
Fix $C>1$. By Theorem~\ref{thm:metricembedding}, there are a Banach space $\XX$, an increasing directed set $(\VV_\lambda)_{\lambda\in \Lambda}$ consisting of finite-dimensional subspaces of $\XX$ such that $\VV:=\cup_{\lambda\in \Lambda} \VV_\lambda$ is dense in $\XX$, and for each $\lambda\in \Lambda$ a linear projection $P_\lambda\colon \XX\to\XX$ with $P_\lambda(\XX)=\VV_\lambda$ and $\norm{P_\lambda}\le C$. For each $0<p\le 1$, the linear map
\[
T_p=\colon\F_p(\XX)\mapsto \YY_p:=\enpar{\oplus_{\lambda\in \Lambda}\F_p(\VV_\lambda)}_{\ell_\infty},
\quad\mu\mapsto (L[P_\lambda;p](\mu))_{\lambda\in \Lambda},
\]
is bounded by $C$. Besides, if $\mu\in\Pt(\Mt)$ is supported on $\VV$, then it is supported of $\VV_\lambda$ for some $\lambda\in\Lambda$, whence
\[
\norm{T_p(\mu)}\ge\norm{L(P_\lambda;p)(\mu)}\ge\norm{L[\VV_\lambda,\XX;p](L(P_\lambda;p)(\mu))}=\norm{\mu}.
\]
Consequently, $T_p$ is an isomorphic embedding by \cite{AACD2018}*{Proposition 4.17}. In fact,
\[
\norm{\mu} \le \norm{T(\mu)} \le C\norm{\mu}
\] for all $\mu\in \F_p(\XX)$. In particular $T_p$ is one-to-one. In turn, combining Theorem~\ref{thm:FpVSchauder} with Proposition\ref{prop:SchauderSeparates} gives that the linear contraction
\[
S_{p}\colon\YY_p\to\YY_1,
\quad
(\mu_\lambda)_{\lambda\in\Lambda}\mapsto (J_{p}[\VV_\lambda](\mu_\lambda))_{\lambda\in\Lambda},
\]
is one-to-one. Finally, $L[\Mt,\XX;p]$ is one-to-one by Theorem~\ref{thm:CR}. Since the diagram
\[
\xymatrixcolsep{5pc}\xymatrix{\F_p(\Mt)\ar[r]^{L[\Mt,\XX;p]}\ar[d]_{J_{p}[\Mt]}
&\F_p(\XX)\ar[r]^{T_p}\ar[d]_{J_{p}[\XX]}
&\YY_p\ar[d]^{S_{p}} \\
\F(\Mt)\ar[r]_{L[\Mt;\XX;1]}
&\F(\XX)\ar[r]_{T_1}
&\YY_1}
\]
commutes, $J_{p}[\Mt]$ is one-to-one.
\end{proof}
\subsection*{Conflict of Interest}
The authors declare that they have no conflict of interest.

\end{document}